\documentclass[12pt]{article}
\usepackage{amssymb}
\newtheorem{theorem}{Theorem}[section]

\usepackage{amsfonts,amssymb,eucal,amsmath,cases}
\title{Some Applications of Generalized Mountain Pass Lemma }

\author{\small\sc Fengying Li\footnote{Email:lify0308@163.com} \\
 \small \it  The School of Economic and Mathematics,
Southwestern
University of Finance and Economics, \\
\small \it Chengdu 611130, China\\
\small\sc \ Bingyu Li and \ Shiqing Zhang \\
 \small \it College of
Mathematics, Sichuan University,
 \small\it Chengdu 610064, People's Republic of China}
\date{}

\begin{document}
\maketitle

\begin{abstract} The Ghoussoub-Preiss's
generalized  Mountain Pass Lemma  with Cerami-Palais-Smale type
condition is a generalization of classical MPL of
Ambrosetti-Rabinowitz, we apply it to study the existence of the
periodic solutions with a given energy for some  second order
Hamiltonian
systems with symmetrical and non-symmetrical potentials.\\
\end{abstract}



\section*{1. Introduction and Main Results}
\setcounter{section}{1} \setcounter{equation}{0}

In 1948, Seifert([17]) studied the periodic solutions of the
Hamiltonian systems using geometrical and topological methods; in
1978 and 1979, Rabinowitz([15,16])studied the periodic solutions of
the Hamiltonian systems using global variational methods; in 1980's,
Benci ([4])and Gluck-Ziller([9]) and Hayashi([11]) used Jacobi
metric and very complicated geodesic methods and algebraic topology
to study the periodic solutions for second order Hamiltonian systems
with a fixed energy:
\begin{numcases}{}
\ddot{q}+V'(q)=0\label{1.1}\\
\frac{1}{2}|\dot{q}|^2+V(q)=h\label{1.2}
\end{numcases}

They proved  the following very general theorem:
\begin{theorem}\label{th:1.2}
\:Suppose $V\in C^1(R^n,R)$ ,if
$$\{x\in R^n|V(x)\leq h\}$$
 is bounded, and
 $$V^{\prime}(x)\not=0,\ \ \ \ \forall x\in\{x\in R^n|V(x)=h\},$$
then the (1.1)-(1.2) has a periodic solution with energy h.
\end{theorem}

   For the existence of multiple periodic solutions for (1.1)-(1.2),
 we can refer Groessen([10]) and Long [12] and the references there.\\

   Ambrosetti--Coti Zelati([1]) used Ljusternik-Schnirelmann theory with classical $(PS)^+$ compact
condition to get the following Theorem:
\begin{theorem}\label{th:1.2}
\:Suppose $V \in C^2(\mathbb R^n \backslash \{0\},\mathbb R)$
satisfies:

$(A1)$.\:$3V'(x)\cdot{x}+V''(x)x\cdot{x}\neq0,\:\forall\,x\in\Omega=
\mathbb{R}^n\backslash\{0\}$;

$(A2)$.\:$V'(x)\cdot{x}>0,\quad\forall\,x\in\Omega$;

$(A3')$.\:$\exists\,\alpha\in(0,2)$, such that
$$
V'(x)\cdot{x}\geq-\alpha V(x),\quad \forall\:x\in\Omega;
$$

$(A4')$.\:$\exists\,\delta\in(0,2)$ and $r>0$,\:such that
$$
V'(x)\cdot{x}\leq-\delta V(x),\quad \forall\:0 <|x|\leq r;
$$

$(A5')$.\:$\underset{|x|\rightarrow+\infty}{\liminf}\left[V(x)
+\dfrac{1}{2}V'(x)\cdot{x}\right]\geq 0$.

Then $\forall\,h<0$ the system (1.1)-(1.2) has at least a
non-constant weak periodic solution which satisfies (1.1)-(1.2)
pointwise except on a zero-measurable set.
\end{theorem}
  Ambrosetti-Coti Zelati ([2]) used a variant of the classical Mountain-Pass Lemma
and a constraint minimizing method to get the following Theorems:
\begin{theorem} \label{th:1.3}
Suppose $V\in
C^1(\mathbb{R}^n\backslash\{0\},\mathbb{R})$ satisfies: \\

$(V1)$.\:$V(-\xi)=V(\xi),\:\forall\,\xi\in\Omega=\mathbb{R}^n\backslash\{0\}$;

$(V2)$.\:$\exists\,\alpha\in[1,2)$,\:such that
$$
\nabla V(\xi)\cdot{\xi}\geq -\alpha V(\xi)>0,\quad
\forall\,\xi\in\Omega;
$$

$(V3)$.\:$\exists\,\delta\in(0,2)$ and $r>0$,\:such that
$$
\nabla V(\xi)\cdot{\xi}\leq - \delta V(\xi),\quad
\forall\:0<|\xi|\leq r;
$$

$(V4)$.\:$V(\xi)\rightarrow 0$,\qquad as \: $|\xi|\rightarrow +\infty$.\\

Then $\forall\,h<0$, the problem $(1.1)-(1.2)$ has a weak periodic solution.
\end{theorem}
\begin{theorem} \label{th:1.4}
Suppose $V$ satisfies  $(V1),(V3),(V4)$ and

$(V2')$.\:$\exists\,\alpha \in(0,2)$,\: such that
$$
\nabla V(\xi)\cdot{\xi}\geq - \alpha V(\xi)>0,\quad
\forall\,\xi\in\Omega;
$$

$(V5)$.\:$V\in C^2(\Omega,\mathbb{R})$ and
$$
3\nabla V(\xi)\cdot{\xi} + V''(\xi)\xi\cdot{\xi}>0.
$$

Then $\forall h<0,(1.1)-(1.2)$ has a weak periodic solution.
\end{theorem}

Yuan-Zhang([19]) proved the following Theorem:
\begin{theorem}\label{th:1.5}
Suppose $V\in C^1(\mathbb R^n\backslash \{0\},\mathbb R)$ satisfies:

$(V_1)$.\: $V(-q)=V(q);$

$(V_2)$.\:There are constant $0<\alpha<2 $ such that
$$
\langle V'(q),q\rangle \geqslant -\alpha V(q)>0,\quad \forall\, q\in
\mathbb {R}^n\backslash  \{0\};
$$

$(V_3)$.\: $\exists\, \delta\in (0,2), r>0$, such that
$$
\langle V'(q),q\rangle\leqslant-\delta V(q),\quad \forall\,0<|q|\leq
r;
$$

$(V_4)$. $V(q)\rightarrow0$,\quad as\:\: $|q|\rightarrow +\infty$.

Then for any given $h<0$, the system (1.1)-(1.2) has at least a
non-constant weak periodic solution which can be obtained by
Mountain Pass Lemma.
\end{theorem}

Motivated by these papers ,we use Ghoussoub-Preiss's Generalized Mountain Pass Lemma with
Cerami-Palais-Smale condition at some levels for a closed subset to
 study the new periodic solutions with symmetrical and
non-symmetrical potentials, we obtain the following Theorems:

{\bf Theorem 1.6}\ \ Suppose $V\in C^1(R^n,R)$ and $h\in R$
satisfies

$(B_1)$\ \ $ V(-q)=V(q).$

$(B_2)$\ \ $\exists \mu_1>0,\mu_2\geq 0,s.t.$ $V^{\prime}(q)\cdot
q\geq\mu_1V(q)-\mu_2.$

$(B_3)$\ \ $V(q)\geq h, |q|\rightarrow +\infty.$

$(B_4)$\ \ $\forall q\not=0,3V^{\prime}(q)\cdot q+V''(q)q\cdot q\not=0.$

Then for any $h>\frac{\mu_2}{\mu_1},$ $(1.1)-(1.2)$ has at least
one non-constant periodic solution with the given energy h, which can be obtained by
the generalized MPL method.\\

{\bf Corollary1.1}\ \ Suppose $a>0,\mu_1\geq 2,\mu_2\geq 0,V(q)=a|q|^{\mu_1}+\frac{\mu_2}{\mu_1}$,
then the conditions of Theorem1.1 hold and for any $h>\frac{\mu_2}{\mu_1}$ ,
$(1.1)-(1.2)$ has at least
two non-constant periodic solution with the given energy h.\\

{\bf Theorem 1.7}\ \ Suppose $V\in C^1(R^n,R)$ and $h\in R$
satisfies ($B_2$), ($B_3$) and
($B_5$)\ \ $\exists r>0$, s.t.
$$\inf_{u\in F}\int_0^1(h-V(u))dt>0,$$
where $$F\triangleq\{u\in H^1|\|\dot{u}\|_{L^2}=r\}.$$ Then $\forall
h>\frac{\mu_2}{\mu_1}$, (1.1)-(1.2) has at least one non-constant
periodic solution with energy $h$.

\section{ A Few  Lemmas}
\setcounter{section}{2} \setcounter{equation}{0}

 Define Sobolev space:
 $$H^1=W^{1,2}(R/TZ,R^n)=\{u:R\rightarrow R^n,u\in L^2,\dot{u}\in L^2,u(t+1)=u(t)\}$$
 Then the standard $H^1$ norm is equivalent to
 $$\|u\|=\|u\|_{H^1}=\left(\int^1_0|\dot{u}|^2dt\right)^{1/2}+|\int_0^1 u(t)dt|.$$

{\bf Lemma 2.1}([1,10])\ \ Let
$f(u)=\frac{1}{2}\int^1_0|\dot{u}|^2dt\int^1_0(h-V(u))dt$ and
$\widetilde{u}\in H^1$ be such that $f^{\prime}(\widetilde{u})=0$
 and $f(\widetilde{u})>0$. Set

\begin{equation}
\frac{1}{T^2}=\frac{\int^1_0(h-V(\widetilde{u}))dt}{\frac{1}{2}\int^1_0|\dot{\widetilde{u}}|^2dt}\label{2.1}
\end{equation}
Then $\widetilde{q}(t)=\widetilde{u}(t/T)$ is a non-constant
$T$-periodic solution for (1.1)-(1.2).

By symmetry condition $(B_1)$, similar to Ambrosetti-Coti
Zelati[2], let
$$E_1=\{u\in
H^1=W^{1,2}(R/Z,R^n),u(t+1/2)=-u(t)\},$$
$$E_2=\{u\in
H^1=W^{1,2}(R/Z,R^n),u(-t)=-u(t)\}.$$

 By the symmetrical condition $(B_1)$ and Palais's symmetrical principle([14])
 or similar proof of [1,2],we have\\

{\bf Lemma 2.2}\ \ If $\bar{u}\in E_i$ is a critical point of
$f(u)$ and $f(\bar{u})>0$, then $\bar{q}(t)=\bar{u}(t/T)$ is a
non-constant $T$-periodic solution of
(1.1)-(1.2).\\

 Using the famous Ekeland's variational principle, Ekeland proved

{\bf Lemma 2.3}(Ekeland[7])\ \ Let $X$ be a  Banach space, $F\subset
X$ be a closed (weakly closed) subset.
 Suppose that $\Phi $ defined on $F$ is Gateaux-differentiable and lower
 semi-continuous (or weakly lower semi-continuous)
  and bounded from
 below. Then there is a sequence $x_n\subset F$ such that
 $$\Phi(x_n)\rightarrow\inf_{F}\Phi$$
 $$(1+\|x_n\|)\|\Phi^{'}(x_n)\|\rightarrow 0.$$

Motivated by the paper of Cerami[6], Ekeland [7],
Ghoussoub-Preiss[8] presented a weaker compact condition than the
classical $(CPS)_c$ condition:

 {\bf Definition 2.1}([7,8])\ \ Let $X$ be a  Banach
space, $F\subset X$ be a closed  subset, let $\delta(x,F)$ denotes
the distance of $x$ to the set $F$.
 Suppose that $\Phi $ defined on $X$ is Gateaux-differentiable,
 if sequence $\{x_n\}\subset X$ such that
 $$\delta(x_n,F)\rightarrow 0,$$
 $$\Phi(x_n)\rightarrow c,$$
 $$(1+\|x_n\|)\|\Phi^{'}(x_n)\|\rightarrow 0,$$

then $\{x_n\}$ has a strongly convergent subsequence.

Then we call $f$ satisfies $(CPS)_{c,F}$ condition at the level
$c$ for the closed subset $F\subset X$,
we denote it as $(CPS)_{c,F}$\\

We can give a weaker condition than $(CPS)_c$ condition:

{\bf Definition 2.2}\ \ Let $X$ be a  Banach space. $F\subset X$ be
a weakly closed subset.
 Suppose that $\Phi $ defined on $X$ is Gateaux-differentiable,
 if sequence $x_n$ such that
$$\delta(x_n,F)\rightarrow 0,$$
 $$\Phi(x_n)\rightarrow \gamma,$$
 $$(1+\|x_n\|)\|\Phi^{'}(x_n)\|\rightarrow 0,$$

then $\{x_n\}$ has a weakly convergent subsequence.

Then we call $f$ satisfies $(WCPS)_{c,F}$ condition.\\

Now by $\bf Lemma2.3$, it's easy to prove

{\bf Lemma 2.4}\ \ Let $X$ be a  Banach space,

(i). Let $F\subset X$ be a closed subset. Suppose that $\Phi $
defined on $X$ is Gateaux-differentiable and lower semi-continuous
and bounded from below, if $\Phi$ satisfies $(CPS)_{\inf\Phi,F}$
condition, then $\Phi$ attains its infimum on $F$.

 (ii). Let $F\subset X$ be a weakly closed subset. Suppose that $\Phi $ defined on $F$ is Gateaux-differentiable and weakly lower semi-continuous
 and bounded from below,
 if $\Phi$ satisfies $(WCPS)_{inf\Phi,F}$ condition, then $\Phi$ attains its infimum on $F$.

 {\bf Definition 2.3([7,8])}\ \ Let $X$ be a  Banach space,
 $F\subset X$ be a closed  subset. If $z_0,z_1$ belong different
 disjoint connected components in $X\backslash F$, then we call $F$
 separates $z_0$ and $z_1$.

 Motivated by the famous classical Mountain Pass Lemma of
 Ambrosetti-Rabinowitz [3], Ghoussoub-Preiss[8] gave a generalized
 MPL:

{\bf Lemma 2.5 (Ghoussoub-Preiss's generalized MPL [8],[7])}\ \ Let
$X$ be a Banach space.Suppose that $\Phi(u):X\rightarrow R$  is  a
continuous Gateaux-differentiable function with
$\Phi^{\prime}:X\rightarrow X^{*}$ norm-to-weak$^{*}$ continuous.
Take two points $z_0,z_1$ in $X$, and define
$$\Gamma=\{c\in C^0([0,1];X)|c(0)=z_0, c(1)=z_1\}$$
$$\gamma=\inf_{c\in\Gamma}\max_{0\leq t\leq 1}\Phi(c(t))$$
Let $F\subset X$ be a closed subset separating $z_0$ and $z_1$.
Assume that

$$ \Phi(x)>\max\{\Phi(z_0),\Phi(z_1)\},\forall x\in F,$$
$\Phi$ satisfies condition $(CPS)_{\gamma,F}$ on the level $\gamma$
for the set $F$. Then there is a critical point of $\Phi$ on the
level $\gamma.$

\section{The Proof of Theorem 1.6}
\setcounter{section}{3} \setcounter{equation}{0}

We define weakly closed subsets of $H^1$:
$$F=\{u\in H^1|\int_0^1(V(u)+\frac{1}{2}V^{\prime}(u)u)dt=h\}.$$
$$F_i=\{u\in E_i|\int_0^1(V(u)+\frac{1}{2}V^{\prime}(u)u)dt=h\},i=1,2.$$

{\bf Lemma 3.1} If $(B_2)-(B_4)$ hold,then $F,F_1,F_2\not=\emptyset$.

{\bf Proof}\ \ Similar to the proof of [1].Let $u\in H^1,u\not=0$ be fixed.
For $a>0$,let\\
$$g_u(a)=g(au)=\int_0^1[V(au)+\frac{1}{2}V'(au)au]dt$$

By $(B_4)$,$\frac{d}{da}g_u(a)\not=0$,so $g_u$ is  strictly monotone.
Notice that\\
$$g_u(0)=g(0)=V(0)\leq\frac{\mu_2}{\mu_1}$$\\
When $a$ is large,we use $(B_2)-(B_3)$ to have

$$g_u(a)=g(au)=\int_0^1[V(au)+\frac{1}{2}V'(au)au]dt$$\\
$$\geq(1+\frac{\mu_1}{2})\int_0^1V(au)dt-\frac{\mu_1}{2}$$\\
$$\geq(1+\frac{\mu_1}{2})h-\frac{\mu_1}{2}$$\\
Hence $\forall h>\frac{\mu_2}{\mu_1}$, we have
$$g_u(+\infty)=g(+\infty)>h$$
So for any given $u\in H^1,u\not=0$,there is $a(u)>0$ such that $a(u)u\in F$.
Similarly we can prove that for any given $u\in E_i,u\not=0$,there is $a(u)>0$ such that $a(u)u\in F_i.$

 {\bf Lemma 3.2} If $(B_1),(B_2)$ and $(B_4)$ hold , then for any given $c>0$,
 $f(u)$ satisfies $(CPS)_{c,F_i}$
 condition, that is :
 If $\{u_n\}\subset H^{1}$ satisfies
\begin{equation}
\delta(u_n,F_i)\rightarrow 0,
 f(u_n)\rightarrow c>0,\ \ \ \
(1+\|u_n\|)f^{\prime}(u_n)\rightarrow 0.\label{3.1}
\end{equation}
Then $\{u_n\}$ has a strongly convergent subsequence.\\

{\bf Proof}\ \ Notice that $\forall u\in E_i, \int_0^1 u(t)dt=0$, so
we know $\|u\|_{E_i}\triangleq(\int_0^1|\dot{u}|^2dt)^{1/2}$ is an
equivalent norm on $E_i$. Now from $f(u_n)\rightarrow c$, we have
\begin{equation}
-\frac{1}{2}\|u_n\|_{E_i}^2\cdot\int^1_0V(u_n)dt\rightarrow
c-\frac{h}{2}\|u_n\|_{E_i}^2\label{3.2}
\end{equation}

By $(B_2)$ we have

\begin{eqnarray}
<f^{\prime}(u_n),u_n>&=&\|u_n\|_{E_i}^2\cdot\int^1_0(h-V(u_n)-\frac{1}{2}<V^{\prime}(u_n),u_n>)dt\nonumber\\
&\leq&\|u_n\|_{E_i}^2\int^1_0[h+\frac{\mu_2}{2}-(1+\frac{\mu_1}{2})V(u_n)]dt\label{3.3}
\end{eqnarray}

By (3.2) and (3.3) we have
\begin{eqnarray}
<f^{\prime}(u_n),u_n>&\leq&(h+\frac{\mu_2}{2})\|u_n\|_{E_i}^2+(1+\frac{\mu_1}{2})(2c-h\|u_n\|_{E_i}^2)\nonumber\\
&=&(-\frac{\mu_1}{2}h+\frac{\mu_2}{2})\|u_n\|_{E_i}^2+C_1\label{3.4}
\end{eqnarray}
Where $C_1=2(1+\frac{\mu_1}{2})c$

Since $h>\frac{\mu_2}{\mu_1}$, then (3.1)and (3.4) imply $\|u_n\|_{E_i}$
is bounded.

The rest for proving $\{u_n\}$ has a strongly convergent subsequence
is standard.

{\bf Remark 3.1} We notice that in our proof, we didn't use the
condition
\begin{equation}
\delta(u_n,F_i)\rightarrow 0.
\end{equation}
 It seems interesting to efficiently
use this condition to weak our assumptions.

{\bf Lemma 3.2} Let\\
\begin{equation}
G=\{u\in H^1|\int_0^1(V(u)+\frac{1}{2}V^{\prime}(u)u)dt<h\},
\end{equation}
\begin{equation}
G_i=\{u\in E_i|\int_0^1(V(u)+\frac{1}{2}V^{\prime}(u)u)dt<h\}.
\end{equation}
Then

(i).$F,F_i,i=1,2$ are respectively the boundaries of $G,G_i$.

(ii).If $(B_1)$ holds, then $F,F_i,G,G_i$ are symmetric with respect to
the origin $0$.

(iii).If $V(0)<h$ holds, then $0\in G,G_i,i=1,2.$\\

It's not difficult to prove the following two Lemmas:

 {\bf Lemma 3.3} $f(u)$ is weakly lower semi-continuous on
 $H^1$ and $F,F_i.$

{\bf Lemma 3.4} $F,F_i,i=1,2.$ are weakly closed subsets in $H^1$.

 {\bf Lemma 3.5} The functional $f(u)$ has
positive lower bound on $F_i.$

{\bf Proof}\  By the definitions of $f(u)$ and $F_i$, we have
\begin{equation}
f(u)=\frac{1}{4}\int^1_0|\dot{u}|^2dt\int^1_0(V^{\prime}(u)u)dt,u\in
F_i.
\end{equation}
For $u\in F_i$ and $(B_2)$ ,we have
$$\frac{1}{2}V'(u)u=h-V(u)\geq h-\frac{1}{\mu_1}V'(u)u-\frac{\mu_2}{\mu_1},$$
$$V'(u)u\geq\frac{h-\frac{\mu_2}{\mu_1}}{\frac{1}{2}+\frac{1}{\mu_1}}>0.$$

  So we have the functional $f(u)\geq 0$.
  Furthermore, we claims that
\begin{equation}
 \inf f(u)>0,
\end{equation}

since otherwise, $u(t)=const$ attains the infimum 0.

 If $u\in F_i$, then  by the
symmetry $u(t+1/2)=-u(t)$ or $u(-t)=-u(t)$, we know $u(t)=0,\forall t$; by ($B_2$)
we have $V(0)\leq\frac{\mu_2}{\mu_1}$, by $h>\frac{\mu_2}{\mu_1}$, we have
$V(0)<h$.
By the definition of $F_i$, $0\notin F_i$. So
\begin{equation}
\inf_{F_2} f(u)>0.
\end{equation}

Now by  Lemmas 3.1-3.5 and Lemma 2.4, we know $f(u)$ attains the
infimum on $F_i$, and we know that the minimizer is nonconstant .\\

{\bf Lemma3.6} $\exists z_1 \not=0,z_1\in H^1$ s.t. $f(z_1)\leq 0.$

 {\bf Proof}\  For any given $y_1\not= const$,$\dot{y}_1\not=0$,so
 $min|y_1(t)|>0,$
 we let $z_1(t)=Ry_1(t)$, then when R is large enough,
 by condition $(B_3)$, we have

\begin{equation}
\int_0^1(h-V(z_1))dt\leq0 ,\label{3.24}
\end{equation}
that is,
\begin{equation}
f(z_1)\leq 0.\label{3.25}
\end{equation}

 {\bf Lemma3.7} $f(0)=0.$

{\bf Lemma3.8} $F_i$ separates $z_1$ and $0$.

 {\bf Proof} By $V(0)<h$, we have that $0\in G_i$.
By $(B_2)$ and $(B_3)$ and $h>\frac{\mu_2}{\mu_1}$, we can choose R large enough
such that
 $$z_1=Ry_1\in \{u\in H^1|\int_0^1(V(u)+\frac{1}{2}V^{\prime}(u)u)dt$$
 $$\geq(1+\frac{\mu_1}{2})\int_0^1V(u)dt-\frac{\mu_1}{2}$$
$$\geq(1+\frac{\mu_1}{2})h-\frac{\mu_1}{2}>h\}$$.

So $F_i$ separates $z_1$ and $0$.

Now by Lemmas 2.4-2.5, 3.1-3.8, we can prove Theorem 1.6.

\section{The Proof of Theorem 1.7}

Let
$$F=\{u\in H^1|\|\dot{u}\|_{L^2}=r\},$$
$$G_1=\{u\in H^1|\|\dot{u}\|_{L^2}<r\},$$
$$G_2=\{u\in H^1|\|\dot{u}\|_{L^2}>r\}.$$
Then $H^1\setminus F=G_1\cup G_2$.

Notice that we can use $(B_5)$ to get that

$$F\cap\{u\in H^1|f(u)\geq c\}=\{u\in H^1,
\frac{1}{2}r^2\int_0^1(h-V(u))dt\geq c\},$$
\begin{eqnarray*}
&H^1&\setminus(F\cap\{u\in H^1|f(u)\geq c\}) \\
&=&\{u\in H^1|\|\dot{u}\|_{L^2}<r\}\cup \{u\in
H^1|\|\dot{u}\|_{L^2}>r\}\cup\{u\in H^1|f(u)<c\}.
\end{eqnarray*}
It's easy to see $u_1=0\in G_1$, we choose $u_2$ such that
$\|\dot{u_2}\|_{L^2}>r$, so $u_2\in G_2$. Now every path $g(t)$
connecting $u_1$ and $u_2$ must pass $F$, so we have
 $$\max_{0\leq t\leq 1}f(g(t))\geq \inf_{u\in F}f(u)=(\frac{1}{2}r^2)\inf_{u\in F}\int_0^1(h-V(u))\geq c>0.$$

 So from the above, in order to apply Ghoussoub-Preiss's generalized
 MPL, now we only need to prove the closed set $F$ separate $u_1$
 and $u_2$ and $f$ satisfies $(CPS)_{c,F}$.

 From the definitions of the set $F$ and $u_1$ and $u_2$, we know
 $F$ separate $u_1$ and $u_2$.

 In order to prove $f$ satisfies $(CPS)_{c,F}$ for any $c>0$,
 firstly, from $(B_2)$, similar to the proof of Lemma 3.1, we can
 get $(\int_0^1|\dot{u}_n|^2dt)^{1/2}$ is bounded, then
 by $(B_3)$, we prove that $|u_n(0)|$ is bounded. In fact,
if otherwise, there exists a subsequence, we still denote it as
$\{u_n(0)\}$ satisfying
$$|u_n(0)|\rightarrow +\infty.$$
 By Newton-Leibniz formula
and Cauchy-Schwarz inequality, we have
\begin{eqnarray}
\min_{0\leq t\leq 1}|u_n(t)|&\geq&|u_n(0)|-\|\dot{u}_n\|_2\rightarrow +\infty\nonumber\\
\label{3.8}
\end{eqnarray}
So by $(B_3)$ we have
\begin{equation}
\int^1_0V(u_n)dt\geq h,\ \ \ \ \rm{as}\ n\rightarrow
+\infty,\label{3.9}
\end{equation}
\begin{equation}
\lim\limits_{n\rightarrow\infty}f(u_n)=\lim\limits_{n\rightarrow\infty}\frac{1}{2}\int^1_0|\dot{u}_n|^2dt\int^1_0(h-
V(u_n))dt\leq 0,\label{3.10}
\end{equation}

 which contradicts with
$f(u_n)\rightarrow c>0$.

 We know that $H^1$ is a
reflexive Banach space, so $\{u_n\}$ has a weakly convergent
subsequence. The rest that proving $\{u_n\}$ has a strongly
convergent subsequence is standard, we can refer to Ambrosetti-Coti
Zelati [2].

\section*{Acknowledgements}

 We would like to thank the supports of NSF of China and a research fund for
 the Doctoral program of higher education of China.

\end{document}